\documentclass[a4paper,11pt]{amsart}

\usepackage{amsmath}
\usepackage{amsfonts}
\usepackage{amssymb}
\usepackage{amsthm}
\usepackage[english]{babel}
\usepackage[latin1,ansinew]{inputenc}
\usepackage{a4wide}
\usepackage{dsfont}

\usepackage{graphicx}
\newtheoremstyle{note}% <name> 
 	{3pt}% <Space above>
 	{3pt}% <Space below> 
 	{\itshape}% <Body font>
 	{}% <Indent amounti>
 	{\bfseries}
 	{.}% <Punctuation after theorem head> 
 	{.5em}% <Space after theorem head>
 	{\thmname{#1}\thmnumber{ \textup{#2}}\thmnote{ \upshape#3}}
\theoremstyle{note}
\newtheorem{thm}{Theorem}%[section]
\newtheorem*{lemma*}{Lemma}

\newtheorem*{proposition*}{Proposition}
\newtheorem*{thm*}{Theorem}
\newtheorem{defi}{Definition}
\newcommand{\R}{{\mathbb R}}

\newcommand{\C}{{\mathbb C}}

\newcommand{\nn}{\nonumber}
\newcommand{\bea}{\begin{eqnarray}}
\newcommand{\eea}{\end{eqnarray}}
\newcommand{\beann}{\begin{eqnarray*}}
\newcommand{\eeann}{\end{eqnarray*}}
\newcommand{\ba}{\begin{array}}
\newcommand{\ea}{\end{array}}

\newcommand{\beq}{\begin{equation}}
\newcommand{\eeq}{\end{equation}}
\newcommand{\be}{\begin{equation}}
\newcommand{\ee}{\end{equation}}

\DeclareMathOperator{\Real}{Re}
\providecommand{\norm}[1]{\lVert#1\rVert}

 \DeclareMathOperator{\sgn}{sgn} 
\date{March 6, 2013}
\title{An analytical proof for the stability of 
Heimburg-Jackson pulses}

\author{Heinrich Freist\"uhler}
\address{Universit\"at Konstanz, 78457 Konstanz, Germany}
\email{heinrich.freistuehler@uni-konstanz.de}
\author{Johannes H\"owing}
\address{Universit\"at Hamburg, 20146 Hamburg, Germany}
\email{johannes.hoewing@math.uni-hamburg.de}


\begin{document}

\begin{abstract}
This paper studies analytically the stability 
of solitary waves 
in a generalized Boussinesq equation with quadratic-cubic nonlinearity.
For general values of two parameters $a$ and  $b$ determining the system, unstable
waves may occur. If however, as in a situation for which this Boussinesq equation was recently 
proposed as a model for pulse propagation in nerves, $(a,b)$ belongs to a certain natural regime,
then all possible waves are stable. 
\end{abstract}
\maketitle
\section{Situation and results}
This note is directly prompted by the
article \cite{HJ} in which Heimburg and Jackson suggest
the partial differential equation
\be\label{BQext}
v_{tt} + (-v + av^2 +b v^3)_{xx} + v_{xxxx} = 0
\ee
%with 
%\be\label{BQext_flux}
%p(V) =\quad\text{with } a,b\in\R.
%\ee
%\eqref{BQext} with \eqref{BQext_flux}
as a model for pulse propagation in biomembranes and nerves and argue 
that this model reflects certain properties of nerve axons better than the 
well known Hodgkin-Huxley and FitzHugh-Nagumo equations. Since the appearance of \cite{HJ}
this model has been studied intensely; for general aspects of these studies we refer the 
reader to the recent survey \cite{ARH}.
The interest in equation \eqref{BQext}
rests on the fact that it admits solitary waves, i.\ e., traveling-wave solutions
\be
v(x,t)=V(x-ct)\quad\text{with }\ V(\pm\infty)=0.
\label{tw}\ee
It is some of these solitary waves that Heimburg and Jackson propose as good representations for pulses 
in the abovementioned biological contexts. 
Now, as in order for this to be the case, the solitary waves should be dynamically stable, 
they and collaborators recently studied this issue computationally \cite{LAJH} 
and found that solitary waves are numerically stable in the case 
that the two parameters $a$ and $b$ occurring in \eqref{BQext} assume certain values 
%$a=8.3$ and $b=-26.5$; these values 
that are significant for the concrete contexts they investigate. \\

The present note gives a complete picture of the existence and stability of solitary waves 
in the {\it extended Boussinesq}\footnote{We use this name in analogy with common terminology for
the extended Korteweg-de Vries equation $v_t+(av^2+bv^3)_x+v_{xxx}=0$ (cf., e.\ g., \cite{GPPS}).}
equation
\eqref{BQext} by analytical deduction.
While the extreme cases $a=0,$ and $b=0$ have been well understood before 
(cf.\ \cite{BS}), no simple scaling 
argument applies to the case $ab\neq 0$.
In fact, the literature does not seem to provide any concrete results concerning 
the stability of solitary waves for {\it generalized 
Boussinesq}\footnote{We use this name in analogy with common terminology for
the generalized Korteweg-de Vries equation $v_t+(p(v))_x+v_{xxx}=0$ (cf., e.\ g., \cite{GSS}).}
 equations
\be\label{BQgen}
v_{tt} + (p(v))_{xx} + v_{xxxx} = 0
\ee
with $p''$ non-monomial.
Ours here rely on findings reported in  
\cite{H}. As \cite{H}, our argumentation follows 
Grillakis, Shatah, Strauss \cite{GSS}
and Bona and Sachs \cite{BS} in  
considering 
the so-called \emph{moment of instability}'s second derivative, 
the sign of which  
% $m(c)$ and its second derivative $m''(c)$ --  
%scalar functions of the wavespeed $c$ which 
allows to conclude or preclude the existence of growing modes in the 
linearization of \eqref{BQext} around a solitary wave \eqref{tw}.
\par\medskip

We first characterize the set of all solitary waves that are possible for equation
\eqref{BQext}.

\begin{thm}\label{T_Existence}
Equation \eqref{BQext} admits positive solitary waves %of elevation with maximum $v_{\max}(c)$ 
of speed $c$ if and only if $a,b,c$ satisfy 
\begin{itemize}
\item $b>0,$ $a\in\R$ and $c^2\in [0,1),$ or  %and $v_{\max}(c):= \bar v(c),$
\item $b<0,$ $a>0$ and $c^2\in \left[\max\left\{0,1+\frac{2a^2}{9b}\right\},1\right).$ %and $v_{\max}(c):=v^+,$
\end{itemize}
It admits negative solitary wave %with minimum $v_{\min}(c)$ 
of speed $c$ if and only if $a,b,c,$ satisfy 
\begin{itemize}
\item $b>0,$ $a\in\R$ and $c^2\in [0,1),$ or %and $v_{\min}(c):=v^-$
\item $b<0,$ $a<0$ and $c^2\in \left[\max\left\{0,1+\frac{2a^2}{9b}\right\},1\right).$ %and $v_{\min}(c):=v^-.$
\end{itemize}
With 
$$
\bar v:=\frac{2}{b}\left(-\frac{a}{3}+\sqrt{\frac{a^2}{9}+\frac{b}{2}(1-c^2)}\right),
\quad\hbox{and}\quad
\underline{v}:=\frac{2}{b}\left(-\frac{a}{3}-\sqrt{\frac{a^2}{9}+\frac{b}{2}(1-c^2)}\right),
$$
each 
positive solitary wave has the respective value $\bar v$ as its maximum, and each negative 
solitary wave has $\underline{v}$ as its minimum. 
\end{thm}
To give a precise definition of stability for this context, 
we write \eqref{BQext} as a system of first order 
in time:
\be\begin{aligned}\label{BQ}
v_t-u_x&=0,\\
u_t+p(v)_x&=-v_{xxx}.
\end{aligned}
\ee
%
% --------------------------Definition Stability-------------------------------
%
\begin{defi}\label{definition1_B} \cite{BS}
A traveling wave $(V,U)$ of \eqref{BQ} is called \emph{(orbitally) stable} if for each 
$\varepsilon>0,$ there exists a $\delta>0$ such that for any solution
$$(v,u)\in (V,U)+C([0,T);H^3(\R)\times H^2(\R))$$ of \eqref{BQ}, 
closeness at initial time, 
$$
\norm{(v,u)(\cdot,0)-(V,U)(\cdot)}_{H^1\times L^2}<\delta
$$
implies, besides existence for all times (i.\ e., one may take $T=\infty$), in particular
orbital closeness at any time,
$$
\inf_{\sigma\in\R}\norm{(v,u)(\cdot,t)-(V,U)(\cdot+\sigma)}_{H^1\times L^2}<\varepsilon
\quad\text{for all }t>0.
\label{foranyt}
$$
\end{defi}

\begin{defi}
We call solitary waves of \eqref{BQext} {\rm Heimburg-Jackson pulses}, if 
\be
b\le-\frac13a^2.
\label{HJineq}
\ee
\end{defi}

The following is the main result of this paper.

\begin{thm}\label{T_HeimburgJackson}
All Heimburg-Jackson pulses are stable. 
\end{thm}

While there is no equivalence, for arbitrary generalized 
Boussinesq equations \eqref{BQgen}, between stability of constant states 
and stability of solitary waves (cf.\ \cite{H}, assertion (ii) of Theorem 4a),
the following seems enlightening for the family of equations under study.

\begin{thm}
Equation \eqref{BQext} is linearly wellposed at any constant state if and only if
\eqref{HJineq} holds.
\end{thm}

In other words, for \eqref{BQext}, stability of constant states does imply stability of all 
solitary waves.
 
\par\medskip
We also show 
\begin{thm}
(i) Assume that $a>0$ and  
\be
b>-\frac29a^2.
\label{twoninths}
\ee
Then there are values $0<c_*\le c^*<1$  such that while all positive waves of speeds with $c^2>c^{*2}$
are stable, all positive waves of speeds with $c^2<c_*^2$ are unstable. Furthermore there are values 
$0<c_\flat\leq c_\sharp<1$ such that all negative waves of  speeds with $c^2<c_\flat^2$ and
all negative waves of  speeds with $c^2>c_{\sharp}^2$ are unstable. \\
(ii) Interchanging the roles of positive and negative waves, the same statement holds given \eqref{twoninths} and $a<0$.
\end{thm}

The transition, for positive waves, between stability for 'fast' waves and instability for 'slow' waves  
vaguely reminds of such a transition in the FitzHugh-Nagumo model, cf.\ \cite{J,KSS}.
\par\medskip
Theorems 1 and 4 imply in particular that for certain choices of $a$ and $b$ violating \eqref{HJineq},
there are unstable solitary waves. \\

Theorems 1, 2, 3, 4 will be demonstrated in Section 2.\\

The following finding is useful for deciding (in-)stability 
of individual solitary waves for cases violating \eqref{HJineq}.

\begin{thm}
For any solitary wave in \eqref{BQext} there is a simple algebraic expression 
$$
\mu(a,b,c)
$$
depending only on the system parameters $a$ and $b$ and the wave's speed $c$ such  
that the wave is stable [unstable] if $\mu(a,b,c)$ is positive [negative].
\end{thm}

Section 3 comprises a proof of Theorem 5 and plots of $\mu$ that also illustrate Theorem 4.

\section{Proofs of Theorems 1 through 4}
As on the one hand the cases $a=0$ and $b=0$ are covered in the literature as mentioned above
and on the other hand the transformation $v\to -v$ is equivalent to replacing $p(v)$ with
$-p(-v)$, we assume without loss of generality for the 
remainder of this paper that 
$$
a>0
\quad\text{and}\quad
b\neq 0.
$$
\subsection*{Proof of Theorem 1} 
\begin{figure}
\includegraphics[scale=0.6]{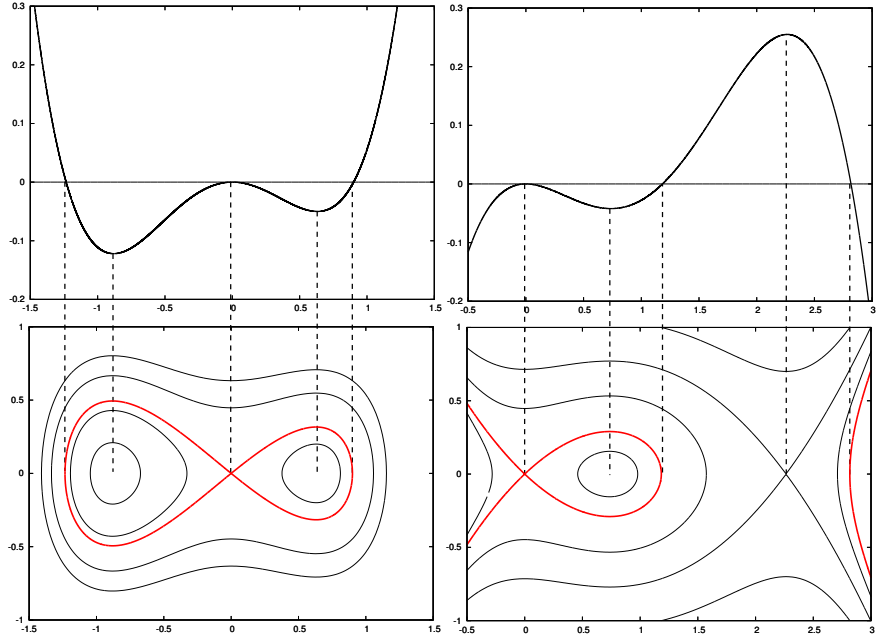}
\caption{Graph of $F$ and level curves of $I$ with $p(v)=-v+av^2+bv^3$ for fixed speed $c=\frac{2}{3}$ for $a=\frac{1}{4},$ $b=1$ (left), and $a=2,$ $b=-\frac{1}{3}$ (right).}
\label{figure1}
\end{figure}

With 
\be
F(v,c)=\frac{1}{2}(c^2-1)v^2+\frac{a}3v^3+\frac{b}4v^4,
%\quad
%-\frac{df(v)}{dv}=p(v),
\ee
a solitary wave satisfies the profile equation
\be\label{profil_B}
\begin{aligned}
V''&=-c^2V - p(V)\\
&=-\frac{\partial F(V,c)}{\partial v};
\end{aligned}
\ee
this equation admits the first integral
\be\nn
\begin{aligned}
I(V,V')&=\frac{1}{2}V'^2 + F(V,c).
%\\
%&= \frac{1}{2}v'^2 +\frac{1}{2}c^2v^2 - f(v)
\end{aligned}
\ee
%with $f$ defined (up to a constant) by
%\be\nn
%-\frac{df(v)}{dv}=p(v).
%\ee
In order for a solitary wave to be at least possible, $(V,V')=(0,0)$ must be a saddle point;
this is the case if and only if $$c^2<1,$$
which we henceforth assume.
Theorem 1 follows directly (cf. Figure \ref{figure1}) from the fact that besides at $0$, $F(.,c)$ vanishes
exactly at $\bar v$ and $\underline{v}$.
\par\medskip

Solitary waves thus occur in families $^cV$ parametrized by their speed $c$.
The key tool for stability considerations is the 
so-called moment of instability,
\be\nn\begin{aligned}
m(c)&=\int_{-\infty}^{\infty}(^cV')^2\;dx,
\end{aligned}
\ee
and our proofs of Theorems 2, 4, and 5 are based on the following fact.

\begin{lemma*}
The solitary wave $^cV$ is stable [unstable] if and only 
if the second derivative $$m''(c)$$ of the moment at the respective speed $c$ is positive [negative].
\end{lemma*}
For this fact and the underlying theory, we refer the reader to \cite{GSS,BS,Z,H}.

\subsection*{Proof of Theorem 2}

Heimburg-Jackson pulses (with $a>0$) are positive. 
As in \cite{H}, we obtain
\be\nn\begin{aligned}
m(c)&=\int_{-\infty}^{\infty}(^cV')^2\;dx=2\int_{0}^{\bar v(c)}(-2 F(v,c))^{1/2}\;dv\\
&=4\int_0^{\sqrt{\bar v(c)}}(-2 F(\bar v(c)-w^2,c))^{1/2}\;w\;dw,\quad\text{with }w:=(\bar v(c)-v)^{1/2}. 
\end{aligned}
\ee
Differentiating twice yields
% \be\begin{aligned}\begin{split}
% \frac{1}{2}m''(c)&=-\int_{v_*}^{\bar v(c)}\frac{(v-v_*)^2+2c(v-v_*)\bar v'(c)}{(-2 F(v,c))^{1/2}}+ \frac{c(v-v_*)^2\bigl[ F_v(v,c)\bar v'(c)+F_c(v,c)\bigr]}{(-2 F(v,c))^{3/2}}\;dv.
% \end{split}\end{aligned}\ee
\be\label{generalm}
m''(c)=2\int_{0}^{\bar v(c)}\frac{v\bigl(2F(v,c)(v+2c\bar v'(c))-cv(F_v(v,c)\bar v'(c)+F_c(v,c))\bigr)}{(-2F(v,c))^{3/2}}\;dv
\ee
with 
\be\nn
\begin{aligned}
F(v,c) &= \frac{1}{2}(c^2-1)v^2+\frac{a}3v^3+\frac{b}4v^4,\\% \frac{1}{2}c^2v^2 -\frac{1}{2}v^2\pm\frac{1}{3}v^3+\frac{k}{4}v^4,\\
F_v(v,c)&=(c^2-1)v +a v^2 + bv^3 = c^2v +p(v),\\
F_c(v,c)&=cv^2,\\
%\bar v(c) &= \frac{1}{3}\frac{-2 + \sqrt{4-18kc^2+18k}}{k}\\
\bar v'(c)&=-\frac{F_c(\bar v(c),c)}{F_v(\bar v(c),c)}.
%= \frac{-c\bar v(c)^2}{c^2\bar v(c)-\bar v(c)+\bar v(c)^2+k\bar v(c)^3}\\
%&=-\frac{c\bar v(c)}{c^2-1+\bar v(c)+k\bar v(c)^2}<0.
\end{aligned}
\ee
It is not difficult to verify that
positivity of the integrand in \eqref{generalm}  is equivalent to positivity of 
\be\nn
Q(v):= \frac{b}{2}v^3+\frac{2}{3}av^2-v+c\bar v'(c)\left(c^2-1+\frac{1}{3}av\right).
\ee
Now, one easily checks that 
$$
Q(\bar v(c))=0\quad\hbox{and}\quad Q'(v)<-\left(1+\frac{8}{27}\frac{a^2}{b}\right)\ \hbox{for all }v\in\R.
$$
This implies that $Q$ is indeed positive on the interval $(0,\bar v(c))$ and thus that $m''(c)>0$.

\subsection*{Proof of Theorem 3} 
Equation \eqref{BQgen} is linearly wellposed at constant states
$v_0\in\R$ if and only if every solution of the form 
$w(x,t) = \exp\left( \lambda t + i\omega x\right),$ $\omega\in\R, \lambda\in\C,$  of its linearization
\be\nn
w_{tt} + p'(v_0)w_{xx} + w_{xxxx}= 0
\ee
has $\Real \lambda \leq 0.$  
Since for any  such mode $w$, 
\be\nn
\lambda^2 - p'(v_0) \omega^2 + \omega^4 = 0,
\ee
this is characterized by $p'(v_0)\equiv-1+2av_0+3bv_0^2\leq0.$ $p'$ has no real zeros iff 
\eqref{HJineq} holds.
\subsection*{Proof of Theorem 4} Here, we have to consider two different cases. 
Consider first the case of a positive wave. Instability of standing waves and hence 'slow'ly traveling waves follows from the following observation: At $c=0,$  
\be\label{standingwaves}
m''(0) = 4\int_0^{\bar v(0)}\frac{v^2F(v,0)}{(-2F(v,0))^{3/2}}\;dv<0,
\ee
since $F(.,c)<0$ in the intervall $(0,\bar{v}(c)).$ Continuity of integral and integrand implies then stability of waves with speed $c^2\approx0;$ this observation is actually a special case of \cite{H2}. On the other hand, to prove stability of 'fast' waves, we apply  Theorem 4 in \cite{H}; 
translated into the present situation, this theorem guarantees existence of a $c^{*}\in(0,1)$ 
such that all solitary waves of speed $c^2\in(c^{*2},1)$ are stable, provided that $p'(0)<0$ and $p''(0)>0;$ with 
\be\nn
p'(0)=-1\quad\text{and}\quad p''(0)=2a,
\ee
this is obviously satisfied.\\
Consider now the case of a negative solitary wave, i.e., $b>0$ and
$\underline{v}(c)=\min\, ^cV<0.$ 
The considerations for $c^2 \lesssim 1$ slightly change
as the moment of instability  is now
\be\nn
m(c) = \int_{-\infty}^{\infty} (^cV')^2\;dx = 2\int_0^{\bar v(c)}
\left(-2G(v,c)\right)^{1/2}\;dv\quad\hbox{ with now }\bar{v}(c):=-\underline{v}(c)
\ee
and
\be\nn
G(v,c) = F(-v,c) = \frac{1}{2}(c^2-1)v^2-\frac{a}{3}v^3+\frac{b}{4}v^4,
\ee
and its second derivative is \eqref{generalm} with $F$ and its
derivatives replaced by $G$ and its derivatives. %and with now$$. 
An obvious analogue of relation \eqref{standingwaves} keeps
implying instability of waves with speed close to $0$. The
following observation now shows instability of 'fast' waves with speed $c^2\lesssim1.$ The quantities $\bar v=-\underline{v},m, m''$ extend to the limiting value $c^2=1$,
with $\bar{v}(\pm 1)=(4a)/(3b)$ and
$\bar{v}'(\pm1)=-3/a$. Thus, for all $v\in(0,\bar v(\pm1))$,
\be\nn
\begin{aligned}
\sgn m''(\pm1)& =  \sgn v\left(2G(v,1)\left(v+2\bar{v}'(1)\right) -
v\left(G_v(v,1)\bar{v}'(1)+G_c(v,1)\right)\right)\\
&= \sgn\left(\frac{-2a}{3} + \frac{b}{2}v\right)\\
&<0.
\end {aligned}
\ee
Now by continuity of
integral and integrand this implies $m''(c)<0$ for $c^2\lesssim1.$ 
(Note in passing that in this case the minimum of the wave, and thus the wave's amplitude,
do not tend to zero for $c^2\to 1$). \\

\begin{figure}[ht]
 \includegraphics[scale=0.75]{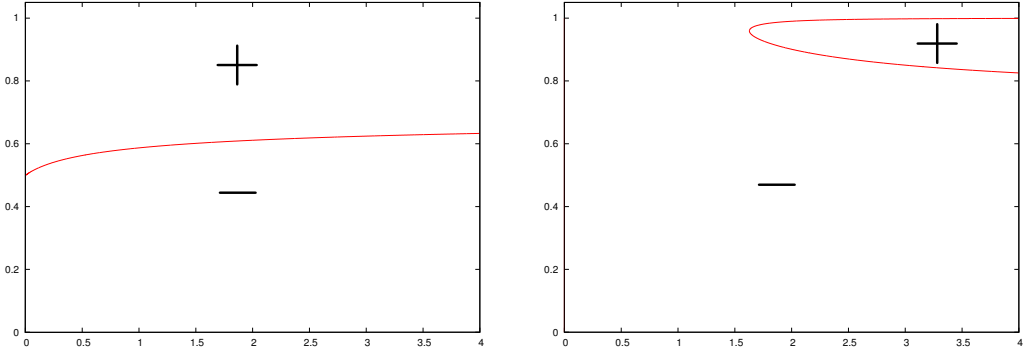}
\caption{Plot of $m''=0$ in the case $a,b>0$ for positive (left) and negative (right) waves. 
The horizontal axis is $k={b}/{a^2}$ and the vertical axis is $c$.
The $\pm$ signs refer to regions where $m''\gtrless0.$}
\label{figure2}
\end{figure}

As can be seen in Figure \ref{figure2}, stable negative waves occur when $k=b/a^2$ is large enough.

\section{Proof of Theorem 5 and plots of \rm{sgn}$(m''(c))$}
We turn from estimating $m''(c)$ to evaluating this quantity.
Let
\be\nn
\begin{aligned}
g(k,c)&:=\arcsin\left(  \frac{2}{3}\frac{1}{\sqrt{\frac{4}{9}+2k(1-c^2)}} \right),\\
\tilde{g}(k,c)&:= -\ln\left(-\frac{\sqrt{\frac{4}{9}+2k(1-c^2)}}{-\frac{2}{3} + \sqrt{-2k(1-c^2)}}\right),
\end{aligned}
\ee  
and 
\be\nn
h(k,c):=\frac{4\left(18kc^2-2-9k\right)\sqrt{1-c^2}}{k\left(2+9k(1-c^2) \right)}.
\ee
\begin{proposition*}
If $a>0$, the following holds with $k=b/a^2$.\\
(i) In the case of $b>0$ and a positive wave, the assertion of Theorem 5 holds with 
\be\nn
\mu(a,b,c) = h(k,c)- 
\frac{4}{3k}\sqrt{\frac{2}{k}}\left(g(k,c)-\frac{\pi}{2}\right).
\ee
(ii) In the case of $b>0$ and a negative wave, the assertion of Theorem 5 holds with 
\be\nn
\mu(a,b,c) = h(k,c)- 
\frac{4}{3k}\sqrt{\frac{2}{k}}\;\left(g(k,c) + \frac{\pi}{2}\right).
\ee
(iii) In the case of $b<0$ and a positive wave, the assertion of Theorem 5 holds with  
\be\nn
\mu(a,b,c) = h(k,c)
-\frac{4}{3k}\sqrt{-\frac{2}{k}}\tilde{g}(k,c).
\ee
\end{proposition*}

{\bf Remark.} Note that with $a>0$, there are no negative solitary waves in the case $b<0$. 
We refrain from formulating the obvious analogue of Proposition 1 for the case $a<0$.  

\begin{proof}
As cases (ii) and (iii) can be treated  analogously, we consider only case (i).
By elementary integration, 
we obtain 
\be\nn
\begin{aligned}
m(c) &= 2\int_0^{\bar{v}(c)} (-2F(v,c))^{1/2}\;dv\\
&= 2\int_0^{\bar{v}(c)} v\left( -\frac{b}{2}v^2-\frac{2}{3}av+1-c^2\right)^{1/2}  \;dv\\
&=\frac{a}{3b^2}\left(2b(1-c^2)+\frac{4}{9}a^2\right)\sqrt{\frac{2}{b}}
%\left(
\arcsin\left( \frac{-bv-\frac{2}{3}a}{\left( \frac{4}{9}a^2+2b(1-c^2) \right)}\right)\Bigg|_{v=0}^{v=\bar v(c)}
%-\arcsin\left(\frac{-\frac{2}{3}a}{\left(\frac{4}{9}a^2+2b(1-c^2)\right)^{1/2}} \right)\right) 
\;+\\
\quad &+ \frac{4}{3b}(1-c^2)^{3/2} + \frac{4a}{9b^2}(1-c^2)^{1/2}.
%\\
%&=\frac{1}{a^2}\left(  \frac{1}{3k^2}\left( 2k(1-c^2)+\frac{4}{9} \right)
%\sqrt{\frac{2}{k}}\;(g(c)-\frac{\pi}2)
% + \frac{4}{3k}\left(1-c^2 \right)^{3/2} +\frac{4}{9k^2}\left(1-c^2\right)^{1/2}\right).
\end{aligned}
\ee
After some slightly tedious calculations, this yields
\be\nn
\begin{aligned}
m''(c) &= 
\frac{4\left(18kc^2-2-9k\right)\sqrt{1-c^2}}{k\left(2+9k(1-c^2) \right)}- 
\frac{4}{3k}\sqrt{\frac{2}{k}}\left(\arcsin\left(  \frac{2}{3}\frac{1}{\sqrt{\frac{4}{9}+2k(1-c^2)}} \right) - \frac{\pi}{2}\right).
\end{aligned}
\ee
\end{proof}

%

%%%%%%%%%%%%%%BIBLIO%%%%%%%%%%%%%%%%%%%%%%

\end{document}